\documentclass[11pt]{amsart}

\usepackage{graphicx}

\newtheorem{proposition}{Proposition}
\newtheorem{theorem}[proposition]{Theorem}

\theoremstyle{definition}
\newtheorem{definition}{Definition}

\def\R{\mathbb{R}}
\def\Z{\mathbb{Z}}
\def\CC{\mathcal{C}}
\def\DD{\mathcal{D}}
\def\J{\mathcal{J}}

\theoremstyle{plain}

\begin{document}

\title{Plane Curves and Contact Geometry}
\author{Lenhard Ng}
\address{Department of Mathematics, Stanford
University, Stanford, CA 94305}
\email{lng@math.stanford.edu}

\begin{abstract}
We apply contact homology to obtain new results in the problem of
distinguishing immersed plane curves without dangerous
self-tangencies.
\end{abstract}

\maketitle

\section{Introduction}

The purpose of this note is to show that contact geometry, and in
particular Legendrian knot theory and contact homology, can be used
to give new information about plane curves without dangerous
self-tangencies. Throughout, the term ``plane curve'' will refer to
an immersion $S^1 \to \R^2$ up to orientation-preserving
reparametrization, i.e., an oriented immersed plane curve in $\R^2$.

\begin{definition}
A self-tangency of a plane curve is \textit{dangerous} if the
orientations on the tangent directions to the curve agree at the
tangency. Two plane curves without dangerous self-tangencies are
\textit{safely homotopic} if they are homotopic through plane curves
without dangerous self-tangencies.
\end{definition}

\begin{figure}
\centerline{
\includegraphics[width=3.8in]{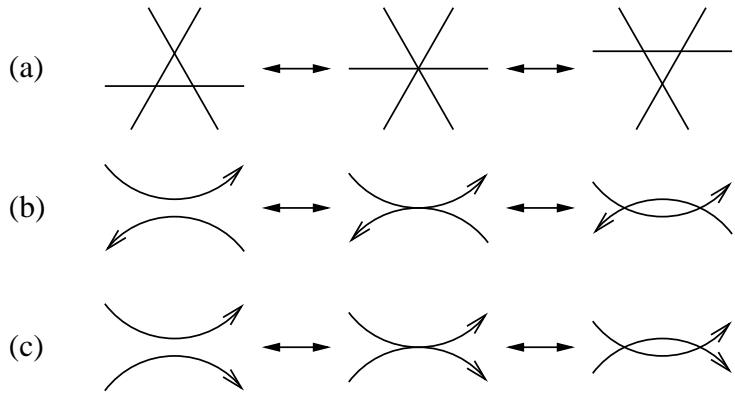}
}
\caption{
Singularities (``perestroikas'') encountered in homotopies of plane
curves: (a) triple point; (b) safe self-tangency; (c) dangerous
self-tangency.
}
\label{fig:perestroikas}
\end{figure}

A generic homotopy of plane curves may contain three types of
singularities, of which one is the dangerous self-tangency; see
Figure~\ref{fig:perestroikas}. Arnold \cite{bib:Ar1,bib:Ar2}
initiated the study of plane curves up to safe homotopy, in
particular introducing a function $J^+$ on plane curves without
dangerous self-tangencies. In the literature, any function of plane
curves without dangerous self-tangencies which does not change under
safe homotopy is called a \textit{$J^+$-type invariant}.

The key point of interest of plane curves without dangerous
self-tangencies is their close link to contact geometry, first noted
by Arnold. There is a natural way to associate to any such plane
curve a Legendrian knot in $\J^1(S^1)$, the $1$-jet space of $S^1$,
which is a contact manifold. We call this the \textit{conormal knot}
of the plane curve. For details, see Section~\ref{ssec:def}.

The conormal knot is a special case of a construction which
associates a Legendrian submanifold to any embedded submanifold of
any manifold, or to any immersed submanifold without dangerous
self-tangencies. This construction has recently been applied to
construct new invariants of knots in $S^3$, and potentially yields
interesting isotopy invariants of arbitrary submanifolds; see
\cite{bib:EE} or \cite{bib:Ng} for an introduction.

There are several well-known $J^+$-type invariants of plane curves,
all arising from the conormal knot construction. The simplest is the
Whitney index, or the degree of the Gauss map of the plane curve.
This is invariant under safe homotopy since it is invariant more
generally under regular homotopy; it also counts the number of times
the conormal knot winds around the base of the solid torus
$\J^1(S^1)$.

A more nontrivial $J^+$-type invariant, as observed by Arnold, is
simply the knot type of the conormal knot in the solid torus. More
interesting still, since the conormal is Legendrian, the contact
planes along the conormal knot give it a framing, and so the framed
knot type of the conormal knot is invariant under safe homotopy. The
framing is measured by a number which is Arnold's original $J^+$
invariant.

To the author's knowledge, all previous work on $J^+$-type
invariants is based on studying the framed knot type of the conormal
knot. For instance, Goryunov \cite{bib:Gor} examined the space of
finite type invariants of plane curves without dangerous
self-tangencies, and Chmutov, Goryunov, and Murakami \cite{bib:CGM}
introduced a $J^+$-type invariant in the form of a HOMFLY polynomial
for the framed conormal knot.

On the other hand, two safely homotopic plane curves have conormal
knots which are isotopic not just as framed knots, but as Legendrian
knots. We will see that the Legendrian type of the conormal knot
gives a finer classification of plane curves than the framed knot
type. The fact (essentially) that Legendrian isotopy is a subtler
notion than framed isotopy was famously demonstrated by Chekanov
\cite{bib:Che} for knots in $\R^3$, using a combinatorial form of
Legendrian contact homology \cite{bib:Eli}. In this paper, we show
that contact homology gives a similar result in our case.

\begin{theorem}[see Propositions~\ref{prop:distinguish}
and~\ref{prop:connectsum}] There are (arbitrarily many) plane curves
with the same framed conormal knot type which are not safely
homotopic.
\label{thm:main}
\end{theorem}

\noindent In the language of Legendrian knot theory, we can rephrase
this result: there are arbitrarily many plane curves whose conormal
knots all have the same classical invariants but are not Legendrian
isotopic.

\begin{figure}
\centerline{
\includegraphics[width=5.3in]{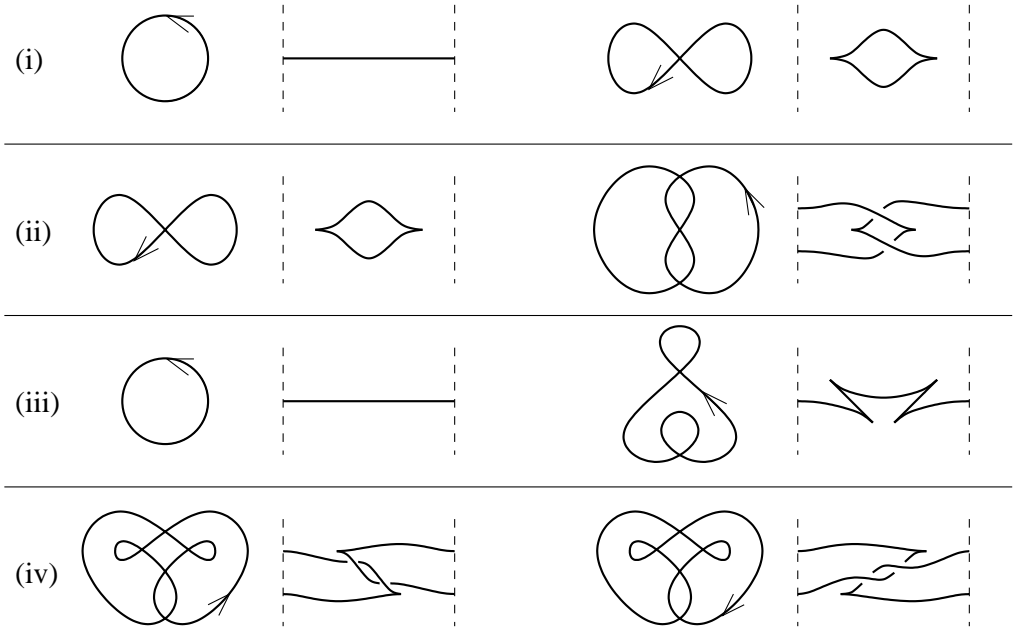}
}
\caption{
Pairs of plane curves, along with their conormal knots, that are
distinguished by increasingly subtle invariants: (i) Whitney index;
(ii) conormal knot type; (iii) framed conormal knot type (Arnold's
$J^+$ invariant); (iv) Legendrian conormal knot type.
}
\label{fig:distinguish}
\end{figure}

An example of a pair of plane curves satisfying the conditions in
Theorem~\ref{thm:main} is given by the bottom line of
Figure~\ref{fig:distinguish}. Note that this ``pair'' is actually
the same plane curve but with different orientations.
Proposition~\ref{prop:distinguish} uses contact homology to
distinguish between these curves.

We review definitions in Section~\ref{ssec:def}, and present an
algorithm for drawing conormal knots in
Section~\ref{ssec:algorithm}. Section~\ref{ssec:distinguish} gives
the proof of our main result, Theorem~\ref{thm:main}. In
Section~\ref{ssec:loops}, we show that contact homology, via recent
work by K\'alm\'an, gives new information about loops of plane
curves as well.

\subsection*{Acknowledgments}

I am grateful to Tobias Ekholm and John Etnyre for useful
discussions, and to Vladimir Chernov for correcting a substantial
error in Section~\ref{ssec:loops}. This work was supported by a
Five-Year Fellowship from the American Institute of Mathematics.

\newpage
\section{Results and Proofs}
\label{sec:conormal}

\subsection{The conormal knot}
\label{ssec:def}

Let $C$ be a plane curve. At each point $x\in C$, the orientation on
$C$ determines two unit vectors, $v_x$ in the direction of $C$ and
$w_x$ given by rotating $v_x$ $90^{\circ}$ counterclockwise.

\begin{definition}
The \textit{conormal knot} of $C$ is the subset of the unit
cotangent bundle $ST^*\R^2$ given by
\[
\{\xi\in ST^*\R^2\,|\,\xi \text{ lies over some } x\in C \text{ and
} \langle \xi,v_x\rangle = 0, \langle \xi,w_x\rangle = 1\}.
\]
The conormal knot inherits an orientation from the orientation on
$C$, since each point on $C$ yields one point in the conormal knot.
\end{definition}

Here the metric on the fibers of $T^*\R^2$ used to define $ST^*\R^2$
is dual to the standard metric on $\R^2$. If $C$ has no dangerous
self-tangencies, then its conormal knot is embedded in $ST^*\R^2$,
and so it makes sense to use the term ``knot.'' We remark that the
conormal knot is actually one half of the usual unit conormal bundle
over the plane curve; the orientation of the plane curve, along with
the orientation of $\R^2$, induces a coorientation on the curve,
which picks out half of the conormal bundle.

The space $ST^*\R^2$ has a natural contact structure given by the
kernel of the $1$-form $\alpha = p_1 \, dq_1 + p_2 \, dq_2$, where
$q_1,q_2$ are coordinates on $\R^2$ and $p_1,p_2$ are dual
coordinates in the cotangent fibers. It is easy to check that the
conormal knot $K$ of any plane curve is Legendrian with respect to
this contact structure, i.e., that $\alpha|_K=0$.

Topologically, $ST^*\R^2 \cong S^1\times\R^2$ is a solid torus, and
it will be more useful for us to view it as the $1$-jet space
$\J^1(S^1) \cong T^*S^1\times\R$. If we set coordinates $\theta,y,z$
on $\J^1(S^1) \cong (\R/2\pi\Z)\times\R\times\R$, then $\J^1(S^1)$
has a natural contact form $\alpha = dz-y\,d\theta$. We can identify
$ST^*\R^2$ and $\J^1(S^1)$ by setting $\theta = \arg(p_1+ip_2)$ (the
argument of the vector $(p_1,p_2)$), $z =
q_1\cos\theta+q_2\sin\theta$, $y=-q_1\sin\theta+q_2\cos\theta$; this
map identifies the contact structures as well.

It is convenient to picture a Legendrian knot in $\J^1(S^1)$ in
terms of its \textit{front}, or projection to $(\R/2\pi\Z)\times\R$
given by the $\theta z$ coordinates. In the subject, ``front'' is
sometimes used in a different sense, namely as a cooriented plane
curve with cusps; for clarity, we will avoid this connotation. A
generic Legendrian knot has a front whose only singularities are
double points and cusps. We can recover a Legendrian knot from its
front by setting $y=dz/d\theta$; in particular, there is no need to
specify over- and undercrossing information for a front. We depict
$(\R/2\pi\Z)\times\R$ by letting $\theta$ be the horizontal axis and
$z$ the vertical axis, and drawing dashed vertical lines to
represent the identified lines $\theta=0$ and $\theta=2\pi$. See
Figure~\ref{fig:distinguish} for examples of fronts in $\J^1(S^1)$.

Any front in $\J^1(S^1)$ has three ``classical'' invariants under
Legendrian isotopy. The first is the knot type of the front in
$\J^1(S^1)$, obtained by smoothing cusps and resolving crossings in
the usual way. The other two are the Thurston--Bennequin number $tb$
and rotation number $r$:
\begin{gather*}
tb = \# \,
\raisebox{-0.13in}{\includegraphics[width=0.3in]{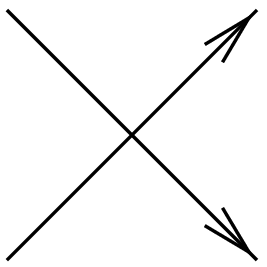}} \, + \,
\#
\,\raisebox{0.17in}{\includegraphics[width=0.3in,angle=180]{tb1.eps}}
\, - \, \#
\,\raisebox{-0.13in}{\includegraphics[width=0.3in,angle=90]{tb1.eps}}
\, - \, \#
\,\raisebox{0.17in}{\includegraphics[width=0.3in,angle=270]{tb1.eps}}
\, - \,
\# \, \raisebox{-0.13in}{\includegraphics[width=0.3in]{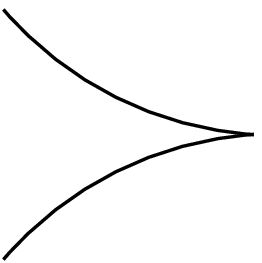}} \\
r = \frac{1}{2} \left( \, \# \,
\raisebox{0.17in}{\reflectbox{\includegraphics[width=0.3in,angle=180]{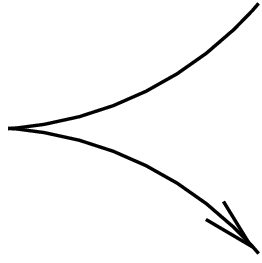}}}
\, + \, \# \,
\raisebox{0.17in}{\includegraphics[width=0.3in,angle=180]{tb3.eps}}
\, - \, \# \,
\raisebox{-0.13in}{\includegraphics[width=0.3in]{tb3.eps}} \, - \,
\# \,
\raisebox{-0.13in}{\reflectbox{\includegraphics[width=0.3in]{tb3.eps}}}
\, \right) .
\end{gather*}
We note that the Thurston--Bennequin number in $\J^1(S^1)$ was first
introduced by Tabachnikov \cite{bib:Tab}.

We now examine the front $K$ of the conormal knot of a plane curve
$C$. There is a simple description for $K$: any point $(q_1,q_2)$ on
$C$, with unit tangent vector $(\cos\varphi,\sin\varphi)$, gives the
point $(\theta,z) = (\varphi+\pi/2,-q_1\sin\varphi+q_2\cos\varphi)$
in $K$, and $K$ is obtained by allowing $(q_1,q_2)$ to range over
$C$. Points of inflection of $C$ correspond to cusps of $K$, and it
is easy to check that any right cusp of $K$ is traversed upwards and
any left cusp downwards; just draw a neighborhood of an inflection
point of $C$.

As for the classical Legendrian invariants of $K$, since $K$ has
equal numbers of left and right cusps, it follows that $r(K)=0$. The
Thurston--Bennequin number of $K$ measures framing and is
essentially Arnold's $J^+$ invariant: $tb(K) = J^+(K)+n(K)^2-1$,
where $n(K)$ is the winding number of $K$ around $S^1$. Hence the
framed knot type of $K$ determines all classical information about
$K$.

\subsection{Drawing the conormal knot front}
\label{ssec:algorithm}

We have already discussed how to define the conormal knot front of a
plane curve, but the definition is not very useful computationally.
Here we present an algorithm for easily obtaining a front isotopic
to the conormal knot front.

Call a plane curve \textit{rectilinear} if it is completely composed
of line segments parallel to either coordinate axis, along with
arbitrarily small smoothing $90^{\circ}$ corners, and no two line
segments lie on the same (horizontal or vertical) line. Clearly any
plane curve is isotopic to a rectilinear curve, and so it suffices
to describe the conormal front for any rectilinear curve.

\begin{figure}
\centerline{
\includegraphics[width=4.5in, bb=148 418 533 663]{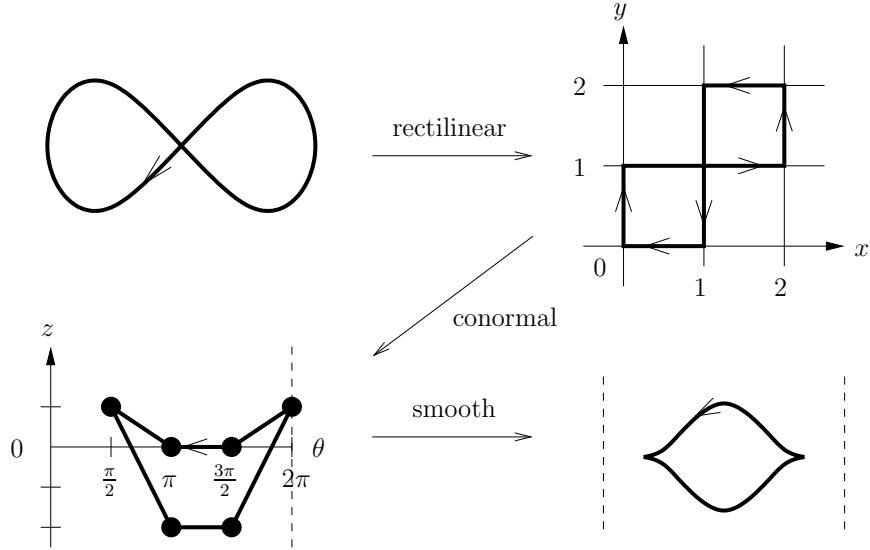}
}
\caption{
Algorithm for obtaining the conormal knot from a plane curve.
}
\label{fig:figure8}
\end{figure}

For ease of notation, label the coordinate axes $x$ and $y$ rather
than $q_1$ and $q_2$. To each line segment $L$ in a rectilinear
plane curve, we associate the following point in
$(\R/2\pi\Z)\times\R$:
\begin{itemize}
\item
$(\pi/2,y)$ if $L$ is in the $+x$ direction and $y$ is the $y$
coordinate of $L$;
\item
$(\pi,-x)$ if $L$ is in the $+y$ direction and $x$ is the $x$
coordinate of $L$;
\item
$(3\pi/2,-y)$ if $L$ is in the $-x$ direction and $y$ is the $y$
coordinate of $L$;
\item
$(0,x)$ if $L$ is in the $-y$ direction and $x$ is the $x$
coordinate of $L$.
\end{itemize}
Next, ``connect the dots'' by joining the points corresponding to
line segments which share an endpoint. Finally, smooth the result,
rounding corners and placing cusps where necessary. See
Figure~\ref{fig:figure8}.

\begin{proposition}
The resulting front in $\J^1(S^1)$ is Legendrian isotopic to the
front of the rectilinear plane curve.
\end{proposition}

\begin{proof}
It is clear that the conormal front for the rectilinear curve passes
through the points given in the algorithm above, since they comprise
the conormal for the line segments of the rectilinear curve. The
conormals of the smoothing corners interpolate between these points.
The conormal front for a smoothing corner at the point $(x,y)$ is
given by $\{(\theta,x\cos\theta+y\sin\theta)\}$ for some range of
$\theta$ in an interval of length $\pi/2$. Hence the conormal fronts
for any two smoothing corners intersect either once or not at all.
It follows that, up to Legendrian isotopy, the conormals for the
smoothing corners can be approximated by the line segments joining
points in the algorithm above.
\end{proof}

\subsection{Nonhomotopic plane curves}
\label{ssec:distinguish}

We can use the algorithm from the previous section to show that
there are plane curves whose conormal knots have the same framed
knot type but which are not Legendrian isotopic.

\begin{proposition}
The plane curves in the bottom line of Figure~\ref{fig:distinguish}
have the same framed conormal knot type but are not safely
homotopic.
\label{prop:distinguish}
\end{proposition}

\begin{figure}
\centerline{
\includegraphics[width=6in]{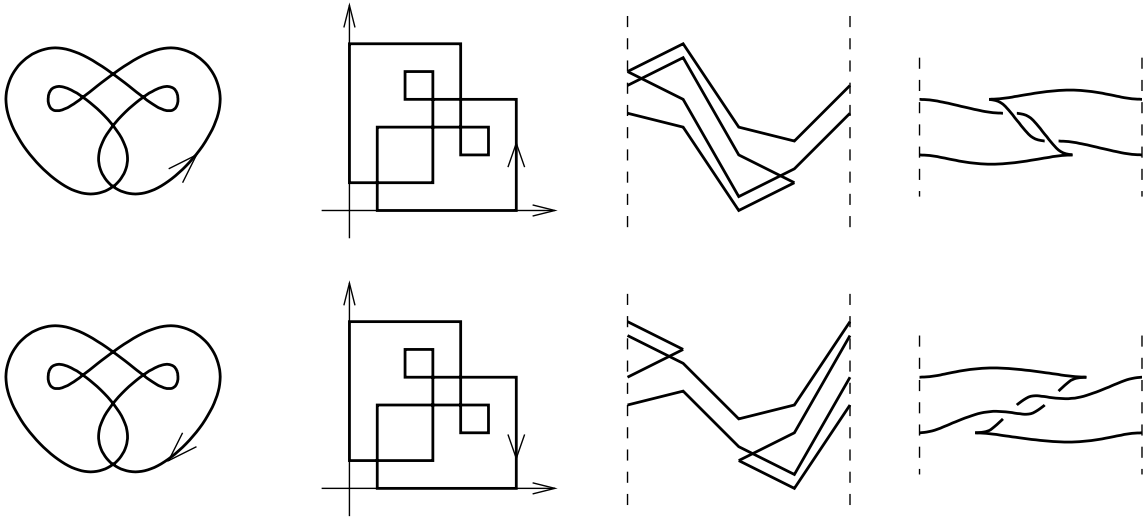}
}
\caption{
Nonhomotopic plane curves, their rectilinear approximations,
conormals, and smoothed conormal fronts.
}
\label{fig:pair}
\end{figure}

\begin{proof}
The two plane curves give conormal knots which are topologically
Whitehead links; see Figure~\ref{fig:pair}. Both conormal knots have
$tb=-3$ (equivalently, $J^+=-2$).

\begin{figure}
\centerline{
\includegraphics[width=5in]{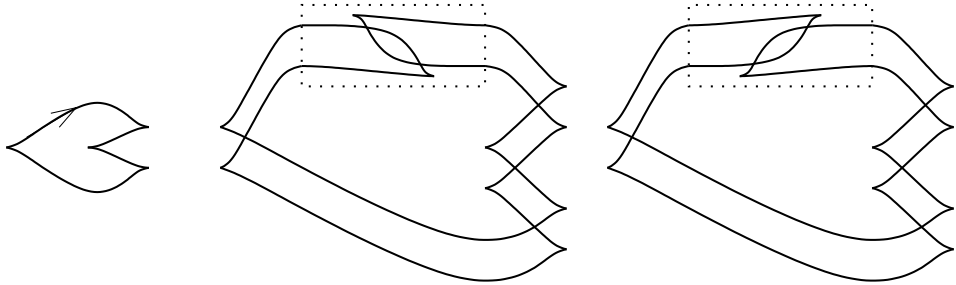}
}
\caption{
The Legendrian satellites of the conormal knots from
Figure~\ref{fig:pair} to the stabilized unknot produce nonisotopic
Legendrian knots.
}
\label{fig:satellite}
\end{figure}

We claim that the two conormal knots are not Legendrian isotopic.
Indeed, applying the Legendrian satellite construction (see the
appendix of \cite{bib:NT}) to the conormal fronts and the stabilized
unknot in $\R^3$ yields two familiar Legendrian knots in $\R^3$:
these are called ``Eliashberg knots'' in \cite{bib:EFM} and labeled
$E(2,3)$ and $E(1,4)$. See Figure~\ref{fig:satellite}. The two knots
can be distinguished by their contact homology differential graded
algebras \cite{bib:Che}; in particular, $E(2,3)$ has Poincar\'e
polynomial $2t+t^{-1}$ and $E(1,4)$ has Poincar\'e polynomial
$t^3+t+t^{-3}$. It follows that the two conormal knots in
$\J^1(S^1)$ are not Legendrian isotopic, as desired.
\end{proof}

\begin{figure}
\centerline{
\includegraphics[width=4in]{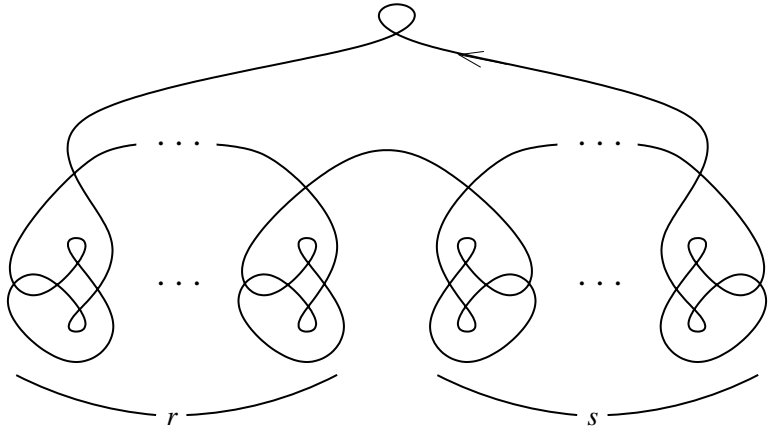}
}
\caption{
The ``connected sum'' plane curve $C_{r,s}$.
}
\label{fig:connectsum}
\end{figure}

We can use the plane curves from Proposition~\ref{prop:distinguish}
to produce an arbitrarily large family of plane curves whose
conormal knots have the same classical invariants but are not
Legendrian isotopic. For $r,s\geq 0$, let $C_{r,s}$ be the plane
curve shown in Figure~\ref{fig:connectsum}, which can be viewed as a
connected sum of the plane curves from
Proposition~\ref{prop:distinguish}. (Note however that the connected
sum operation on plane curves is not well-defined.)

\begin{proposition}
For fixed $n\geq 1$, the $n$ plane curves $C_{r,s}$, $r+s=n$, have
\label{prop:connectsum}
the same framed conormal knot type but are not safely homotopic.
\end{proposition}

\begin{figure}
\centerline{
\includegraphics[width=5in,bb=149 566 596 666]{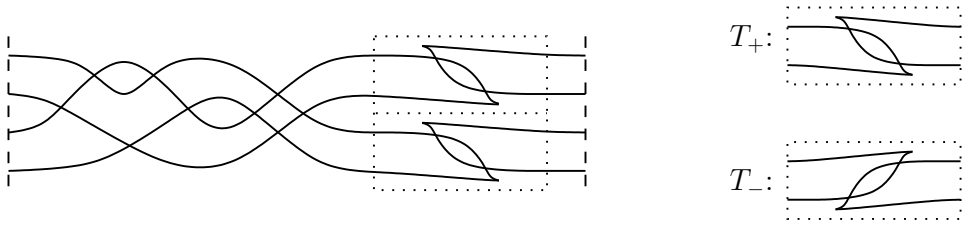}
}
\caption{
The conormal knot for $C_{2,0}$. To obtain the conormal knots for
$C_{1,1}$ and $C_{0,2}$, replace one or both of the boxed tangles
$T_+$ by $T_-$.
}
\label{fig:connectsum-front}
\end{figure}

\begin{proof}[Sketch of proof]
We can use the algorithm from Section~\ref{ssec:def} to find the
conormal fronts for $C_{r,s}$. When $r+s=n$ is fixed, the conormal
fronts for $C_{r,s}$ are identical except for $n$ tangles. Of these
tangles, $r$ are given by the tangle $T_+$ defined in
Figure~\ref{fig:connectsum-front}, and $s$ by $T_-$. The situation
for $n=2$ is shown in Figure~\ref{fig:connectsum-front}; the picture
for $n>2$ is very similar. Note that the conormal fronts for
$C_{r,s}$ are all isotopic as framed knots.

Now consider the Legendrian satellite $K_{r,s}$ of the conormal
front for $C_{r,s}$ to the stabilized unknot, as in the proof of
Proposition~\ref{prop:distinguish}. We distinguish between the knots
$K_{r,s}$ using Poincar\'e polynomials for contact homology.

It is easy to show that $K_{r,s}$ has a graded augmentation, for
instance because it has a ruling \cite{bib:Fu}. An examination of
Maslov indices shows that all crossings in $K_{r,s}$ have degrees
$0,\pm 1,\pm 2$, except for the crossings in the tangles $T_{\pm}$;
the two crossings in any $T_+$ have degree $1$ and $-1$, while the
two crossings in any $T_-$ have degree $3$ and $-3$. Since the
(linearized) differential of the degree $3$ crossing in any $T_-$ is
$0$, we conclude that any Poincar\'e polynomial for $K_{r,s}$ has
$t^3$ coefficient equal to $s$. It follows that the Legendrian knots
$K_{r,s}$, $r+s=n$, are not Legendrian isotopic, and thus that the
plane curves $C_{r,s}$ are not safely homotopic.
\end{proof}

\subsection{Loops of plane curves}
\label{ssec:loops}

Here we consider loops in the space of plane curves. Let $\CC$
denote the space of plane curves, and let $\DD \subset \CC$ be the
discriminant of plane curves with dangerous self-tangencies. We will
present a loop which is contractible in $\CC$ but noncontractible in
$\CC\setminus\DD$.

%
%
%
%
%

\begin{figure}
\centerline{
\includegraphics[width=5in]{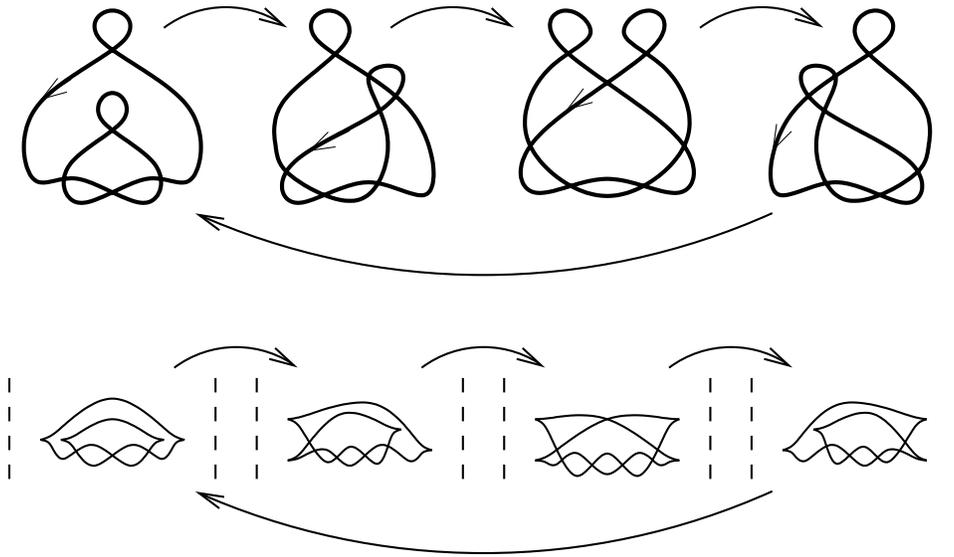}
}
\caption{
A nontrivial loop $\gamma$ of plane curves, and the corresponding
loop $\tilde{\gamma}$ of conormal knots.
}
\label{fig:trefoil-loop}
\end{figure}


Consider the loop $\gamma$ in $\CC\setminus\DD$ pictured in
Figure~\ref{fig:trefoil-loop}. This induces a loop $\tilde{\gamma}$
of Legendrian knots in $\J^1(S^1)$, also shown in
Figure~\ref{fig:trefoil-loop}. As a loop of framed knots in
$\J^1(S^1)$, $\tilde{\gamma}$ essentially rotates the standard
rotationally symmetric diagram of the trefoil by $120^\circ$; since
$\pi_1(SO(3))\cong\Z/2$, $\tilde{\gamma}^6 \sim 1$ in
$\pi_1(\text{framed knots})$. By using the contact condition and
contact homology, we can do even better: a result of K\'alm\'an
\cite{bib:Kal} shows that $\tilde{\gamma}$ has much higher order
when considered as a loop of Legendrian knots.

\begin{proposition}
The loop $\gamma$ is contractible in $\CC$, but has order at least
$30$ in $\pi_1(\CC\setminus\DD)$.
\end{proposition}

\begin{proof}
It is straightforward to check that $\gamma$ is contractible in
$\CC$. Note that by the $h$-principle, $\CC$ is weakly homotopy
equivalent to the space of free loops in $S^1\times\R^2$, which is
not simply connected. However, $\gamma$ can be represented by a loop
of \textit{based} loops in $S^1\times\R^2$, and the space of based
loops is simply connected since $\pi_2(S^1\times\R^2)=0$.

Now consider $\gamma$ as a loop in $\CC\setminus\DD$. The loop
$\tilde{\gamma}$ of Legendrian knots in $\J^1(S^1)$ lifts to an
identical-looking loop $\tilde{\gamma}'$ of Legendrian knots in the
universal cover $\R^3$ with the standard contact structure. (Just
ignore the dashed lines in Figure~\ref{fig:trefoil-loop}.) A result
of \cite{bib:Kal} uses monodromy in contact homology to show that
$\tilde{\gamma}'$ has order at least $30$ in the Legendrian
category; hence $\tilde{\gamma}$ does as well.
\end{proof}


\end{document}